\documentclass[11pt]{article}
\pdfoutput=1

\usepackage[backend=biber,
            bibstyle=phys,
            sorting=none,
            url=false,
            doi=false,
            eprint=true,
            sortcites=true,
            citestyle=numeric-comp]{biblatex}

\DeclareFieldFormat*{title}{\mkbibitalic{#1}}
\DeclareFieldFormat{labelnumberwidth}{\mkbibbrackets{#1}~~}
\renewbibmacro*{volume+number+eid}{%
  \printfield{volume}%
  \iffieldundef{number}{}{%
  \printtext{(}%
  }%
  \printfield{number}%
  \iffieldundef{number}{}{%
  \printtext{)}%
  }%
  \setunit{\addcomma\space}%
  \printfield{eid}}

\usepackage[paperwidth=8.5in, paperheight=11in, portrait, top=1in, bottom=1.in,
 left=1.15in, right=1.15in]{geometry}
\usepackage[utf8]{inputenc}
\usepackage{amssymb,amsthm,amsmath,amsfonts,mathtools,bm,dsfont}
\usepackage{authblk} 
\usepackage{graphicx,mathrsfs,latexsym,xcolor,subfigure,epsfig}
\usepackage{units}
\usepackage[breaklinks]{hyperref}

\newcommand{\qphi}[3]{\varphi_{#1}(#2;#3;q)}
\newcommand{\id}{\mathds{1}}

\def\mU{{\mathcal{U}}}
\def\no{\nonumber}

\numberwithin{equation}{section}
\theoremstyle{plain}
\newtheorem{thm}{Theorem}
\newtheorem{prop}[thm]{Proposition}
\newtheorem{defi}[thm]{Definition}
\newtheorem{rmk}[thm]{Remark}

\newcommand*\pPqskip{8mu}
\catcode`,\active
\newcommand*\pPq{\begingroup
        \catcode`\,\active
        \def ,{\mskip\pPqskip\relax}%
        \dopPq}
\catcode`\,12
\def\dopPq#1#2#3#4#5#6{%
        {}_{#1}\phi_{#2}\biggl(\genfrac..{0pt}{}{#3}{#4};{#5},{#6}\biggr)%
        \endgroup}
\newcommand*\pFqskip{8mu}
\catcode`,\active
\newcommand*\pFq{\begingroup
        \catcode`\,\active
        \def ,{\mskip\pFqskip\relax}%
        \dopFq}
\catcode`\,12
\def\dopFq#1#2#3#4#5{%
        {}_{#1}F_{#2}\biggl(\genfrac..{0pt}{}{#3}{#4};{#5}\biggr)%
        \endgroup}

\title{\bf Bispectrality and biorthogonality of the rational functions of $q$-Hahn type}
\renewcommand*{\Affilfont}{\normalsize\small}
\author[1]{Ismaël Bussière\,}
\author[2]{Julien Gaboriaud\,}
\author[3]{Luc Vinet\,}
\author[4]{Alexei Zhedanov\,\vspace{.5em}}
\affil[1,2,3]{Département de Physique, Université de Montréal, Montréal (Québec), H3C 3J7,
Canada.\vspace{.9em}}
\affil[2,3]{Centre de Recherches Math\'ematiques, Universit\'e de Montr\'eal,
P.O. Box 6128, Centre-ville Station, Montr\'eal (Qu\'ebec), H3C 3J7, Canada.
\vspace{.9em}}
\affil[3]{Institut de valorisation des données (IVADO), Montréal (Québec), H2S 3H1,
Canada.
\vspace{.9em}}
\affil[4]{School of Mathematics, Renmin University of China, Beijing, 100872, China.
\vspace{1.5em}}
{
 \makeatletter
 \renewcommand\AB@affilsepx{: \protect\Affilfont}
 \makeatother
 \affil[ ]{E-mail addresses}
 \makeatletter
 \renewcommand\AB@affilsepx{, \protect\Affilfont}
 \makeatother
 \affil[1]{ismael.olivier.pacifique.bussiere@umontreal.ca}
 \affil[2]{julien.gaboriaud@umontreal.ca}
 \affil[3]{luc.vinet@umontreal.ca}
 \affil[4]{zhedanov@yahoo.com}
}
\date{\today}

\begin{document}
\maketitle

\hrule
\begin{abstract}
We introduce families of rational functions that are biorthogonal with respect to the
$q$-hypergeometric distribution. A triplet of $q$-difference operators $X$, $Y$, $Z$ is
shown to play a role analogous to the pair of bispectral operators of orthogonal
polynomials. The recurrence relation and difference equation take the form of generalized
eigenvalue problems involving the three operators.  The algebra generated by $X$, $Y$, $Z$
is akin to the algebras of Askey--Wilson type in the case of orthogonal polynomials.  The
actions of these operators in three different basis are presented.  Connections with
Wilson's ${}_{10}\phi_9$ biorthogonal rational functions are also discussed.
\\[.5em]
\end{abstract}
\hrule

\section{Introduction}\label{sec:introduction}

The Askey scheme of orthogonal polynomials \cite{KoekoekLeskyetal2010} presents the
classical basic hypergeometric orthogonal polynomials and their $q\to1$ limits in a
hierarchical way, with the most general polynomials (the Askey--Wilson and Wilson
polynomials \cite{AskeyWilson1985}) sitting at the top.  These families of polynomials
possess various properties of interest such as a three term recurrence relation, a
differential/difference equation, a Rodriguez formula, etc.  They also have
interpretations and applications in algebra and group theory as matrix elements and
overlaps in particular.  The algebraic study of their bispectrality through the algebras
of Askey--Wilson type \cite{Zhedanov1991} has shed remarkable light on these orthogonal
polynomials (OP) and led to connections with various fields.

Put simply, the bispectrality of polynomials \cite{Grunbaum2001} corresponds to the fact
that they obey (i) a three term recurrence relation and (ii) a
differential/difference equation that amount to eigenvalue equations for two operators:
(i) one that acts on the degree $n$, with eigenvalue depending on the variable $x$, and
(ii) another one acting on the variable $x$ with eigenvalue involving the degree $n$.

It is quite natural to look for extensions of the $q$-Askey scheme.  One possible
direction is the multivariable  generalization of the various families of the $q$-Askey
scheme (see for example \cite{Tratnik1991, Iliev2011}).  This paper looks at another
direction.  We are interested in biorthogonal rational functions (BRF) that are
bispectral with the aim of enlarging systematically the Askey scheme to BRF.

The biorthogonal rational functions to be considered will be in correspondence
with the classical hypergeometric orthogonal polynomials.  Their biorthogonality relations
will involve the same weights as those for which the polynomials of the Askey scheme are
orthogonal.  The key feature in the case of BRF is that bispectrality will arise through
generalized eigenvalue problems (GEVP) involving three operators.

The simplest case of BRF of the ${}_3F_2$ type was studied in
\cite{TsujimotoVinetetal2020}. These functions which possess the same orthogonality weight
as the Hahn polynomials were called ``BRF of the Hahn type''. Here we shall focus on the
$q$-generalization of the Hahn case. We will obtain and characterize a family of BRF that
are orthogonal with respect to the $q$-hypergeometric distribution and shall call this
family the ``BRF of $q$-Hahn type''.

The study of BRF is not new \cite{RahmanSuslov1993} and various connections with Padé
approximants \cite{Zhedanov2005, SpiridonovZhedanov2007}, continued fractions
\cite{IsmailMasson1995, GuptaMasson1998}, Gram determinants \cite{Wilson1991} and GEVP
\cite{Zhedanov1999, SpiridonovZhedanov2001} have been established before.  Pioneering work
was done by Wilson in \cite{Wilson1991}, where he identified BRF expressed as
${}_{10}\phi_{9}$ basic hypergeometric series.  These functions are very general and
possess a large number of parameters. One can
recover from them the full
$q$-Askey scheme, its $q\to\pm1$ limits and the associated BRF through appropriate limits
and contractions \cite{VandeBultRains2009}.  As will be shown below, both the rational
functions of Hahn and $q$-Hahn type can indeed be obtained from Wilson's BRF through
appropriate limits.

As is well appreciated already, it is most useful to have detailed explicit
characterizations for the various families that descend through limits or specializations
from say the Askey--Wilson polynomials.  This view is adopted here as we study the
rational functions of ${}_{3}\phi_{2}$ type with an eye to working our way up to those of
type ${}_{4}F_{3}$ and ${}_{4}\phi_{3}$.  Let us stress moreover that our focus is on the
bispectral properties as well as their algebraic desription.

Here is the outline for the rest of the paper. In Section \ref{sec:tripletxyz}, the three
operators $X$, $Y$, $Z$ that entail the desired bispectrality will be introduced. The
solutions $U_n$ of the GEVP $YU_n=\lambda_nXU_n$ are derived in Section \ref{sec:gevpfct}
and identified as the BRF of $q$-Hahn type. The biorthgonality of these functions is
the object of Section \ref{sec:ortho} and the precise connection between the BRF of
$q$-Hahn/Hahn type and Wilson's ${}_{10}\phi_{9}$ BRF is detailed in Section
\ref{sec:wilson}.  The tridiagonal actions of $X$, $Y$, $Z$ on the basis of these
functions is presented in Section \ref{sec:tridiag}. How the bispectrality of these
rational functions is encoded in two GEVP is explained in Section \ref{sec:bispectrality}
and the bispectrality algebra realized by the three operators is examined in Section
\ref{sec:algebra}. Final remarks conclude the paper.

\section{A triplet of \textit{q}-difference operators}\label{sec:tripletxyz}
In this paper we assume $q$ to be generic and we use the following notation for
$q$-numbers:
\begin{align}\label{}
 [x]_q=\frac{1-q^x}{1-q}.
\end{align}
It is directly seen that in the limit $q\to1$, the $q$-number $[x]_q$ goes to the
ordinary number $x$.

Let us now introduce the following three $q$-difference operators that will play a central
role in the paper:
\begin{align}
 X^{(\alpha, \beta)}&=[x-\alpha]_q\id-q^{-\alpha}[x]_qT^{-},\label{eq:Xdef}\\
 Y^{(\alpha, \beta)}&=A_1(x)T^++A_2(x)T^-+A_0(x)\id,\label{eq:Ydef}\\
 Z^{(\alpha, \beta)}&=\frac{[-x]_q}{[\alpha-x]_q}T^{-}-\id,\label{eq:Zdef}
\end{align}
where
\begin{align}\label{eq:coeffYshift}
\begin{aligned}
 A_1(x)&=q^{\beta-x}[x-N]_q[x+1-\alpha]_q[x-\alpha]_q,\\
 A_2(x)&=q^{-x}[x]_q[x-N-\alpha+\beta]_q[x-\alpha]_q,\\
 A_0(x)&=-(A_1(x)+A_2(x)).
\end{aligned}
\end{align}
Note that
\begin{align}\label{eq:zfx}
 Z^{(\alpha, \beta)}=-\frac{1}{[x-\alpha]_q}X^{(\alpha, \beta)}
\end{align}
and that $Y^{(\alpha, \beta)}$ has $[x-\alpha]_q$ as a global factor.
We use the notation
\begin{align}
 T^\pm f(x)=f(x\pm 1),\qquad\id f(x)=f(x)
\end{align}
to denote the shift operators $T^{\pm}$ and the identity operator $\id$.

The three operators $X^{(\alpha, \beta)}$, $Y^{(\alpha, \beta)}$, $Z^{(\alpha, \beta)}$
depend on the parameters $\alpha$ and $\beta$. In the remainder of the paper, we will
generally omit this dependency and write $X$, $Y$, $Z$.

These operators act on the space $M$ of real functions $f(x)$ defined on the linear grid
$x=0,\pm 1, \pm 2, \dots$; that is $M=\{f:\mathbb{Z}\to\mathbb{R}\}$. In what follows we
will restrict $M$ to be the space $M_N$ of functions defined on the finite set of points
$x=0,1,\dots,N$. $M_N$ is hence of dimension $N+1$ and one has the following standard
basis on this space:
\begin{align}
 e_k(x)=\delta_{k,x},\qquad k,x\in\{0,1,2,\dots,N\}.
\end{align}
In this basis, the operators $X$, $Y$, $Z$ act as matrices of size $(N+1)\times (N+1)$.
More precisely, the matrices for $X$ and $Z$ are lower bi-diagonal while the matrix for
$Y$ is tridiagonal:
\begin{align}
 X&=\left[
 \begin{matrix}
     [-\alpha]_q&0&0&\dots&0\\
     -q^{-\alpha}[1]_q&[1-\alpha]_q&0&&\vdots\\
     0&\ddots&\ddots&\ddots&0\\
     \vdots&&-q^{-\alpha}[N-1]_q&[N-1-\alpha]_q&0\\
     0&\dots&0&-q^{-\alpha}[N]_q&[N-\alpha]_q
 \end{matrix}
 \right],\\
 Y&=\left[
 \begin{matrix}
     A_0(0)&A_1(0)&0&\dots&0\\
     A_2(1)&A_0(1)&A_1(1)&&\vdots\\
     0&\ddots&\ddots&\ddots&0\\
     \vdots&&A_2(N-1)&A_0(N-1)&A_1(N-1)\\
     0&\dots&0&A_2(N)&A_0(N)
 \end{matrix}
 \right],\\
 Z&=\left[
 \begin{matrix}
    -1&0&0&\dots&0\\
     \tfrac{[-1]_q}{[\alpha-1]_q}&-1&0&&\vdots\\
     0&\ddots&\ddots&\ddots&0\\
     \vdots&&\tfrac{[1-N]_q}{[\alpha-N+1]_q}&-1&0\\
     0&\dots&0&\tfrac{[-N]_q}{[\alpha-N]_q}&-1
 \end{matrix}
 \right].
\end{align}

\section{The rational functions of the ${}_3\phi_2$ type as solutions of a generalized
eigenvalue problem}\label{sec:gevpfct}
In this section we will solve the GEVP built from the pair of operators $X$ and $Y$
\begin{align}
 YU=\lambda XU.
\end{align}
Let us introduce the notation for the basic hypergeometric series ${}_{r+1}\phi_{r}$:
\begin{align}\label{}
 \pPq{r+1}{r}{a_0,a_1,\dots,a_r}{b_1,\dots,b_r}{q}{z}=
 \sum_{k=0}^{\infty}
 \frac{(a_0;q)_k(a_1;q)_k\dots(a_r;q)_k}{(q;q)_k(b_1;q)_k\dots(b_r;q)_k}z^{k},
\end{align}
where $(a;q)_k=(1-a)(1-qa)\dots(1-q^{k-1}a)$ is the $q$-Pochammer symbol.
The main results are summarized in the following statement:
\begin{prop}\label{thm:gevp}
The normalized solutions $U_n$ of the GEVP
\begin{align}\label{eq:gevp}
 YU_n=\lambda_n XU_n
\end{align}
are given by
\begin{align}\label{eq:unorm}
 U_n(x;\alpha,\beta,N;q)=\frac{(q^{-N};q)_n}{(q^{-n-\beta};q)_n}
 \pPq{3}{2}{q^{-n},q^{n+\beta-N},q^{-x}}{q^{-N},q^{\alpha-x}}{q}{q^{\alpha-\beta}},
\end{align}
with the spectrum $\lambda_n$ given by
\begin{align}
 \lambda_n=[-n]_q[n+\beta-N]_q.
\end{align}
The functions $U_n$ will be called ``biorthogonal rational functions of the $q$-Hahn
type''.
\end{prop}
\proof
With hindsight, we introduce the following basis on the space $M_N$ of functions:
\begin{align}\label{eq:qphidef}
 \qphi{n}{x}{\alpha}=\frac{(q^{-x};q)_n}{(q^{\alpha-x};q)_n},\qquad n=0,1,\dots,N.
\end{align}
Enforcing the restriction to $M_N$, it is verified that the operators $X$, $Y$, $Z$ act in
a tridiagonal fashion in this basis:
\begin{align}
 \label{eq:Xphi}
 X\,\qphi{n}{x}{\alpha}&=-[n]_q\qphi{n+1}{x}{\alpha}+[n-\alpha]_q\qphi{n}{x}{\alpha},\\
 \label{eq:Zphi}
 Y\,\qphi{n}{x}{\alpha}&=\nu_n^{(1)}\qphi{n+1}{x}{\alpha}+\nu_n^{(2)}\qphi{n}{x}{\alpha}
  +\nu_n^{(3)}\qphi{n-1}{x}{\alpha},\\
 Z\,\qphi{n}{x}{\alpha}&=\qphi{n+1}{x}{\alpha}-\qphi{n}{x}{\alpha},
 \label{eq:Yphi}
\end{align}
where
\begin{align}
 \nu_n^{(1)}&=-[-n]_q[n]_q[n+\beta-N]_q,\\
 \nu_n^{(2)}&=[-n]_q[n+\beta-N]_q[n-\alpha]_q+q^{\beta-\alpha}[n]_q[n-N-1]_q[1-n]_q,\\
 \nu_n^{(3)}&=-q^{\beta-2\alpha}[n]_q[n-N-1]_q[1+\alpha-n]_q,
\end{align}
and the boundary specializations are given by:
\begin{align}
 X\,\qphi{N}{x}{\alpha}&=[N-\alpha]_q\qphi{N}{x}{\alpha},\\
 Y\,\qphi{0}{x}{\alpha}&=\nu_0^{(1)}\qphi{1}{x}{\alpha}+\nu_0^{(2)}\qphi{0}{x}{\alpha},\\
 Y\,\qphi{N}{x}{\alpha}&=\nu_N^{(2)}\qphi{N}{x}{\alpha}
  +\nu_N^{(3)}\qphi{N-1}{x}{\alpha},\\
 Z\,\qphi{N}{x}{\alpha}&=-\qphi{N}{x}{\alpha}.
\end{align}
Equivalently, the actions of the operators $X$, $Y$, $Z$ on the basis
$\qphi{n}{x}{\alpha}$ correspond to the following matrices:
\begin{align}
 X&=\left[
 \begin{matrix}
     [-\alpha]_q&0&0&\dots&0\\
     0&[1-\alpha]_q&0&&0\\
     0&-[1]_q&[2-\alpha]_q&\ddots&\vdots\\
     &\ddots&\ddots&\ddots&0\\
     \vdots&&-[N-2]_q&[N-1-\alpha]_q&0\\
     0&\dots&0&-[N-1]_q&[N-\alpha]_q
 \end{matrix}
 \right],\\
 Y&=\left[
 \begin{matrix}
     \nu_0^{(2)}&\nu_1^{(3)}&0&\dots&0\\
     \nu_0^{(1)}&\nu_1^{(2)}&\nu_2^{(3)}&&\vdots\\
     0&\ddots&\ddots&\ddots&0\\
     \vdots&&\nu_{N-2}^{(1)}&\nu_{N-1}^{(2)}&\nu_{N}^{(3)}\\
     0&\dots&0&\nu_{N-1}^{(1)}&\nu_{N}^{(2)}
 \end{matrix}
 \right],\\
 Z&=\left[
 \begin{matrix}
     -1&0&0&\dots&0\\
     1&-1&0&&\vdots\\
     0&\ddots&\ddots&\ddots&0\\
     \vdots&&1&-1&0\\
     0&\dots&0&1&-1
 \end{matrix}
 \right].
\end{align}
Remarkably, in the basis of the functions $\qphi{n}{x}{\alpha}$, one readily observes that
the operator $Y$ factorizes as follows:
\begin{align}\label{eq:YXV}
 Y^{(\alpha,\beta)}=X^{(\alpha,\beta)}V^{(\alpha,\beta)},
\end{align}
where $V\equiv V^{(\alpha,\beta)}$ acts as a lower bi-diagonal matrix
\begin{align}
 V\,\qphi{0}{x}{\alpha}&=0,\\
 \label{eq:Vphi}
 V\,\qphi{n}{x}{\alpha}&=[-n]_q[n+\beta-N]_q\qphi{n}{x}{\alpha}
  +q^{\beta-\alpha+1-n}[n]_q[n-N-1]_q\qphi{n-1}{x}{\alpha}.
\end{align}
This factorization leads to a remarkable simplification of the GEVP and allows to find the
an explicit basic hypergeometric solution to the GEVP.  This goes as follows: Expand the
solutions to the GEVP $U_n(x)$ over the set of functions $\qphi{k}{x}{\alpha}$, $k\in
\{0,1,\dots,n\}$. Owing to the factorization \eqref{eq:YXV}, the generalized eigenvalue
equation \eqref{eq:gevp} can be converted into the usual eigenvalue equation
\begin{align}\label{eq:eve}
 VU_n=\lambda_nU_n.
\end{align}
In view of \eqref{eq:Vphi}, the eigenvalues $\lambda_n$ are identified to be
\begin{align}\label{eq:lambdan}
 \lambda_n=[-n]_q[n+\beta-N]_q.
\end{align}
One can now expand $U_n$ in the basis $\qphi{k}{x}{\alpha}$
\begin{align}\label{eq:UnCnk}
 U_n=\sum_{k=0}^{n}C_{n,k}\,\qphi{k}{x}{\alpha}
\end{align}
and it is readily seen using the action \eqref{eq:Vphi} of $V$ on the basis functions
$\qphi{k}{x}{\alpha}$ that \eqref{eq:eve} implies the following two-term
recurrence relation for the coefficients $C_{n,k}$:
\begin{align}
 q^{\beta-\alpha-k}[k+1]_q[k-N]_qC_{n,k+1}
  =\left(-q^{-n}[n]_q[n+\beta-N]_q+q^{-k}[k]_q[k+\beta-N]_q\right)C_{n,k}.
\end{align}
Its solution is found to be
\begin{align}
 C_{n,k}=C_{n,0}\,q^{k(\alpha-\beta)}
  \frac{(q^{-n};q)_k(q^{n+\beta-N};q)_k}{(q;q)_k(q^{-N};q)_k}.
\end{align}
It then follows from \eqref{eq:UnCnk} that
\begin{align}\label{eq:Unx}
 U_n(x;\alpha,\beta,N;q)=\frac{(q^{-N};q)_n}{(q^{-n-\beta};q)_n}
 \pPq{3}{2}{q^{-n},q^{n+\beta-N},q^{-x}}{q^{-N},q^{\alpha-x}}{q}{q^{\alpha-\beta}}.
\end{align}
The normalization constant $C_{n,0}$ has been chosen so that
\begin{align}
\lim_{x\to\infty}U_n=1.
\end{align}
\endproof
\begin{rmk}
It is manifest that the presence of $q$-Pochhammer symbols $\left(q^{\alpha-x};q\right)_n$
in the denominator implies that the function has poles when
$x\in\{\alpha,\alpha+1,\dots,\alpha+n-1\}$.  Hence, one could alternatively express the
function $U_n(x)$ as
\begin{align}
 U_n(x;\alpha,\beta,N;q)=1+\sum_{k=0}^{n-1}\frac{\eta_{n,k}}{[x-\alpha-k]_q}
\end{align}
for some coefficients $\eta_{n,k}$.
\end{rmk}
Note that Askey and Suslov have also studied ${}_{3}\phi_2$ BRF \cite{AskeySuslov1993} but
they do not appear to be related to the ones discussed here.

The functions $U_n(x;\alpha,\beta,N;q)$ will be called ``biorthogonal rational functions
of $q$-Hahn type'' (see Section \ref{sec:ortho} for a justification of the name).
From their hypergeometric expression and contiguity relations, it is seen that the
operators $X$, $Y$, $Z$ shift the parameter $\alpha$ by one:
\begin{align}
 \label{eq:contiguityX}
 X\,U_n(x;\alpha,\beta,N;q)&=[-\alpha]_q\,U_n(x;\alpha+1,\beta,N;q),\\
 \label{eq:contiguityY}
 Y\,U_n(x;\alpha,\beta,N;q)&=[-\alpha]_q[-n]_q[n+\beta-N]_q\,U_n(x;\alpha+1,\beta,N;q),\\
 \label{eq:contiguityZ}
 Z\,U_n(x;\alpha,\beta,N;q)&=-\frac{[-\alpha]_q}{[x-\alpha]_q}\,U_n(x;\alpha+1,\beta,N;q).
\end{align}
Another way to write this is to introduce the operator $P:\alpha\mapsto\alpha+1$.
Then, the operators $X$, $Y$, $Z$ act as
\begin{align}
 X&=[-\alpha]_q\,P,\\
 Y&=[-\alpha]_q[-n]_q[n+\beta-N]_q\,P,\\
 Z&=-\frac{[-\alpha]_q}{[x-\alpha]_q}\,P,
\end{align}
on the solutions $U_n$ of the GEVP.  This will prove useful in the examination of the
biorthogonality properties of these rational functions of Hahn type.

\section{Biorthogonality}\label{sec:ortho}
Let us introduce the following $q$-hypergeometric distribution:
\begin{align}\label{eq:weight}
 w_x^{(\alpha,\beta)}=q^{(1-\alpha)N}\frac{(q^{\alpha-\beta-1};q)_N}{(q^{-\beta};q)_N}
  q^{(\beta+1)x}\frac{(q^{-N};q)_x(q^{1-\alpha};q)_x}{(q;q)_x(q^{\beta-\alpha-N+2})_x},
  \qquad\sum_{x=0}^{N}w_x^{(\alpha,\beta)}=1,
\end{align}
whose normalization is checked with help of the $q$-analog of the Chu--Vandermonde formula
\cite{KoekoekLeskyetal2010}:
\begin{align}\label{}
 \pPq{2}{1}{q^{-N},a}{c}{q}{\frac{cq^{N}}{a}}
  =\frac{(c/a;q)_N}{(c;q)_N}.
\end{align}
This weight is the one with respect to which the $q$-Hahn polynomials are orthogonal
\cite{KoekoekLeskyetal2010}.
\begin{rmk}
Under the involution $S$
\begin{align}\label{}
\begin{aligned}
 S:\alpha&\mapsto\beta-\alpha+2,\\
  x&\mapsto N-x,\\
  q&\mapsto q^{-1},
\end{aligned}
\end{align}
the orthogonality weight is unchanged:
\begin{align}\label{eq:involution}
 S\,w_x^{(\alpha,\beta)}=w_x^{(\alpha,\beta)}.
\end{align}
\end{rmk}
\noindent
We use this $q$-hypergeometric weight to define a scalar product
on the space $\mathcal{L}^{(\alpha,\beta)}$ of rational functions
with parameters $\alpha$ and $\beta$ and defined on the grid $x=0,1,\dots,N$:
\begin{align}\label{eq:scalar}
 \left(f(x),g(x)\right)_{(\alpha,\beta)}
  =\sum_{x=0}^{N}w_x^{(\alpha,\beta)}f(x)g(x),\qquad f,g\in\mathcal{L}^{(\alpha,\beta)}.
\end{align}
In the following, the parameters $(\alpha,\beta)$ will be omitted from the scalar
products.

Our goal is now to find biorthogonal partners for the functions $U_n(x)$, that is,
functions $\mU_m(x)$ orthogonal to $U_n(x)$ under the scalar product \eqref{eq:scalar}:
\begin{align}\label{eq:bio}
 \left(U_n,\mU_m\right)=\delta_{n,m}H_n.
\end{align}
The essential ingredient is the fact that the $U_n$ are solutions to a GEVP.  Given the
scalar product, we can define the adjoint GEVP:
\begin{align}\label{eq:agevp}
 Y^{*}U^{*}_m=\lambda_m X^{*}U^{*}_m,
\end{align}
where the asterisk denotes the adjoint of operators and $U^{*}_m$ denotes the solution of
the adjoint GEVP.

We know that \cite{Zhedanov1999}:
\begin{prop}
The biorthogonal partner of $U_n$ is given by $\mU_m=X^{*}U^{*}_m$.
\end{prop}
\proof
Consider the scalar product $(YU_n,U^{*}_m)$. By definition of the adjoint operators, one
has
\begin{align}
 (YU_n,U^{*}_m)=(U_n,Y^{*}U^{*}_m).
\end{align}
Using the GEVP \eqref{eq:gevp} and the adjoint GEVP \eqref{eq:agevp} leads to
\begin{align}\label{eq:step2}
 \lambda_n(XU_n,U^{*}_m)=\lambda_m(U_n,X^{*}U^{*}_m),
\end{align}
from where we have
\begin{align}
 (\lambda_n-\lambda_m)(U_n,X^{*}U^{*}_m)=0.
\end{align}
In view of formula \eqref{eq:lambdan} for $\lambda_n$, it is clear that first factor is
zero only if $n=m$. Hence the second factor must be zero for $n\neq m$ and this shows the
biorthogonality between $U_n$ and $\mU_m=X^{*}U^{*}_m$.
\endproof

\begin{prop}
Up to a normalization constant $c_{\alpha,\beta}$, the biorthogonal partners
$\mU_m(x;\alpha,\beta,N;q)$ of $U_n(x;\alpha,\beta,N;q)$ are given by
\begin{align}\label{}
 \mU_m(x;\alpha,\beta,N;q)=c_{\alpha,\beta}\, U_m(N-x;\beta-\alpha+2,\beta,N;q^{-1}).
\end{align}
\end{prop}
\proof
We shall obtain $\mU_m=X^*U_m^*$ using the fact that $U_m$ satisfies the GEVP
$YU_m=\lambda_m XU_m$ described in Proposition \ref{thm:gevp} and the properties of the
adjoints of $X$ and $Y$.
The adjoints of the operators involved are computed using the expression for the weight
\eqref{eq:weight} and the definition of the scalar product \eqref{eq:scalar}. A
straightforward calculation yields
\begin{align}
 (T^{-})^{*}&
  =q^{\beta+1}\frac{[x-N]_q[x-\alpha+1]_q}{[x+1]_q[x-N+\beta-\alpha+2]_q}T^{+},\\
 (T^{+})^{*}&=q^{-\beta-1}\frac{[x]_q[x-N+\beta-\alpha+1]_q}{[x-N-1]_q[x-\alpha]_q}T^{-}.
\end{align}
It follows that the adjoints of $X$, $Y$, $Z$ are given by
\begin{align}
 X^{*}&=[x-\alpha]_q\id-q^{\beta-\alpha+1}
  \frac{[x-N]_q[x-\alpha+1]_q}{[x-N+\beta-\alpha+2]_q}T^{+},\\
 Z^{*}&=q^{\beta-\alpha+1}\frac{[x-N]_q}{[x-N+\beta-\alpha+2]_q}T^{+}-\id,\\
 Y^{*}&=q^{-\beta-1}
  \frac{[x]_q[x-N+\beta-\alpha+1]_q}{[x-N-1]_q[x-\alpha]_q}A_1(x-1)T^-\no\\
 &+q^{\beta+1}\frac{[x-N]_q[x-\alpha+1]_q}{[x+1]_q[x-N+\beta-\alpha+2]_q}A_2(x+1)T^+\no\\
 &+A_0(x)\id,
\end{align}
with $A_1(x)$, $A_2(x)$, $A_0(x)$ given by \eqref{eq:coeffYshift}.

We now introduce the transformation $\tau$
\begin{align}\label{}
\begin{aligned}
 \tau:\,&\alpha\mapsto\beta-\alpha+1,\\
  &x\mapsto N-x,\\
  &q\mapsto q^{-1},\\
  &T^{\pm}\mapsto T^{\mp}.
\end{aligned}
\end{align}
to define the operators
\begin{align}\label{eq:tildexy}
 \Tilde{X}=\tau(X),\qquad \Tilde{Y}=\tau(Y).
\end{align}
Note that $\tau$ was chosen in order to reproduce the action of $P\circ S$ on the
parameters $\alpha$, $x$, $q$.

A simple calculation shows that
\begin{align}
 \label{eq:tildex}
 \Tilde{X}\left(\frac{[x-\alpha]_q}{[x-N+\beta-\alpha+1]_q}f(x)\right)
  &=-qX^*\left(f(x)\right),\\
 \label{eq:tildey}
 \Tilde{Y}\left(\frac{[x-\alpha]_q}{[x-N+\beta-\alpha+1]_q}f(x)\right)
  &=-q^{N-\beta+3}Y^*\left(f(x)\right).
\end{align}
Taking the transformation $\tau$ on both sides of \eqref{eq:contiguityX}, we get the
following action of $\Tilde{X}$:
\begin{align}\label{eq:tildexun}
 \Tilde{X}U_n(N-x;\beta-\alpha+1,\beta,N;q^{-1})
  =[\alpha-\beta-1]_qU_n(N-x;\beta-\alpha+2,\beta,N;q^{-1}).
\end{align}
Using \eqref{eq:tildex} and \eqref{eq:tildey}, we can also rewrite the adjoint GEVP
\eqref{eq:agevp} in terms of the transformed operators $\Tilde{X}$ and $\Tilde{Y}$:
\begin{align}\label{}
 -q^{-N+\beta-3}\Tilde{Y}\left(\frac{[x-\alpha]_q}{[x-N+\beta-\alpha+1]_q}U_m^*\right)
  =-q^{-1}\lambda_m\,\Tilde{X}\left(\frac{[x-\alpha]_q}{[x-N+\beta-\alpha+1]_q}
  U^*_m\right),
\end{align}
or equivalently
\begin{align}\label{eq:tildegevp}
 \Tilde{Y}\left(\frac{[x-\alpha]_q}{[x-N+\beta-\alpha+1]_q}U^*_m\right)
  =q^{N-\beta+2}\lambda_m\,\Tilde{X}
  \left(\frac{[x-\alpha]_q}{[x-N+\beta-\alpha+1]_q}U^*_m\right).
\end{align}
We readily find that
\begin{align}
 \tau(\lambda_m)=\tau([-m]_q[m+\beta-N]_q)
  =[-m]_{q^{-1}}[m+\beta-N]_{q^{-1}}
  =q^{N-\beta+2}\lambda_m
\end{align}
and from \eqref{eq:tildegevp} and definition \eqref{eq:tildexy}, we get
\begin{align}\label{}
 \tau(Y)\!\left(\frac{[x-\alpha]_q}{[x-N+\beta-\alpha+1]_q}U^*_m\right)
  =\tau(\lambda_m)\,\,
  \tau(X)\!\left(\frac{[x-\alpha]_q}{[x-N+\beta-\alpha+1]_q}U^*_m\right).
\end{align}
Comparing with the GEVP \eqref{eq:gevp} $YU_m=\lambda_m XU_m$, we see that
\begin{align}\label{}
 \frac{[x-\alpha]_q}{[x-N+\beta-\alpha+1]_q}U^*_m=
 \tau(U_n(x;\alpha,\beta,N;q))=
 U_n(N-x;\beta-\alpha+1,\beta,N;q^{-1}).
\end{align}
We can now obtain the biorthogonal partner $\mU_m=X^*U_m^*$. From this last relation
and \eqref{eq:tildex}, we have
\begin{align*}
 \mathcal{U}_m&=X^*U^*_m \\
 &=-q^{-1}
  \Tilde{X}\left(\frac{[x-\alpha]_q}{[x-N+\beta-\alpha+1]_q}U^*_m\right)\\
 &=-q^{-1}
  \Tilde{X}\left(\frac{[x-\alpha]_q}{[x-N+\beta-\alpha+1]_q}
  \frac{[x-N+\beta-\alpha+1]_q}{[x-\alpha]_q}
  U_m(N-x;\beta-\alpha+1,\beta,N;q^{-1})\right).
\end{align*}
Finally, using \eqref{eq:tildexun}, we get
\begin{align}\label{eq:Vnx}
 \mU_m(x;\alpha,\beta,N;q)
  &=-q^{-1}[\alpha-\beta-1]_q\,U_m(N-x;\beta-\alpha+2,\beta,N;q^{-1}).
\end{align}
\endproof
\noindent
\begin{rmk}
Since $Y$ factorizes as $XV$, the adjoint problem \eqref{eq:agevp} takes the form
$V^{*}X^{*}U^{*}_m=\lambda_mX^{*}U^{*}_m$.
It thus follows that the biorthogonal partners $\mU_m=X^*U_m^*$ satisfy the eigenvalue
equation
\begin{align}
 V^{*}\mathcal{U}_m=\lambda_m\mathcal{U}_m.
\end{align}
\end{rmk}
\noindent
To find the norm $H_n$ in the biorthogonality relation \eqref{eq:bio}, we will make use of
the results in Wilson's paper \cite{Wilson1991} and connect them to ours through an
appropriate limiting procedure. This will be detailed in the next section.

\section{Connection with Wilson's ${}_{10}\phi_{9}$ BRF}\label{sec:wilson}
In his 1991 paper \cite{Wilson1991}, Wilson introduced biorthogonal rational functions
expressed as ${}_{10}\phi_{9}$ basic hypergeometric series.
These functions obey a biorthogonality relation on the $N+1$ points of a discrete grid
$z+z^{-1}$ where $z\in\{q^{a},q^{a+1},\dots,q^{a+N}\}$.
The biorthogonality relation can be written as follows
\footnote{We have corrected two typos from Wilson's formulas: we added a missing $q^{-n}$
and changed $q^{a+1}$ to $q^{2a+1}$ in the expression of $h_n$. These typos were also
observed in \cite{Rosengren2007}.}:
\begin{align}\label{}
 \sum_{x=0}^{N}w_x(a,b,c,d,e,f) u_n(x+a;a,b,c,d,e,f) v_n(x+a;a,b,c,d,e,f)
 = \delta_{nm}h_n(a,b,c,d,e,f),
\end{align}
with
\begin{align}\label{}
 w_x(a,b,c,d,e,f)&=q^{x}\frac{(q^{2a};q)_x}{(q;q)_x}\frac{1-q^{2a+2x}}{1-q^{2a}}\no\\
  &\quad\times\frac{(q^{a+b};q)_x(q^{a+c};q)_x(q^{a+d};q)_x(q^{a+e};q)_x(q^{a+f};q)_x}
  {(q^{1+a-b};q)_x(q^{1+a-c};q)_x(q^{1+a-d};q)_x(q^{1+a-e};q)_x(q^{1+a-f};q)_x},\\
 u_n(x+a;a,b,c,d,e,f)&\no\\
 &\hspace{-7em}=\pPq{10}{9}{q^{-n},q^{\frac{a-e}{2}+1},-q^{\frac{a-e}{2}+1},
   q^{a-e},q^{1-c-e},q^{1-d-e},q^{1-e-b},q^{n-e-f},q^{2a+z},q^{-z}}
 {q^{\frac{a-e}{2}},-q^{\frac{a-e}{2}},q^{a+b},q^{a+c},q^{a+d},q^{1-e+a+n},
   q^{1-n+a+f},q^{1-e-z-a},q^{1-e+z+a}}{q}{q},\\
 v_n(x+a;a,b,c,d,e,f)&=u_n(x+a;b,a,c,d,f,e),\\
 h_n(a,b,c,d,e,f)&=\frac{(q^{2a+1};q)_N(q^{1-c-d};q)_N(q^{1-c-e};q)_N(q^{1-d-e})_N}
 {(q^{b+f};q)_N(q^{1+a-c};q)_N(q^{1+a-d};q)_N(q^{1+a-e};q)_N}\no\\
 &\quad\times q^{-n}(q;q)_n\frac{(q^{n-e-f};q)_n}{(q^{1-e-f};q)_{2n}}
 \frac{(q^{c+d};q)_n(q^{1-e+c};q)_n(q^{1-f+b};q)_n}
 {(q^{a+b};q)_n(q^{-b-e};q)_n(q^{-a-f};q)_n},
\end{align}
with the constraints
\begin{align}\label{eq:constraints}
 q^{a+b+c+d+e+f}=q,\qquad q^{a+b}=q^{-N}.
\end{align}
Note that here we have absorbed normalization factors in the $h_n$ so the functions $u_n$
are no longer symmetric in $a$, $b$, $c$, $d$ whereas it is the case for the functions
$r_n$ in Wilson's paper. Moreover, we have denoted his parameters $z$, $a$, $b$, $c$, $d$,
$e$, $f$ by $q^{x}$, $q^{a}$, $q^{b}$, $q^{c}$, $q^{d}$, $q^{e}$, $q^{f}$ respectively so
as to make the correspondence with our functions more direct.

The limits can now be taken individually for each of the $w_x$, $u_n$, $v_n$, $h_n$.

First, let us account for the constraints \eqref{eq:constraints} and replace $b$ and $f$
using
\begin{subequations}\label{eq:parametrization}
\begin{align}\label{}
 b=-N-a,\qquad f=N+1-c-d-e.
\end{align}
Now, we redefine the remaining parameters as
\begin{align*}\label{}
 a+e=1-\alpha,\qquad c+d=1+\beta
\end{align*}
and replace the parameters $d$ and $e$ in the expressions using
\begin{align}\label{}
 d=1+\beta-c,\qquad e=1-\alpha-a.
\end{align}
\end{subequations}
If we now take the limit $q^{a}\to\infty$, we recover the biorthogonality relation for the
BRF of $q$-Hahn type.
Indeed, under the parametrization \eqref{eq:parametrization}, we obtain
\begin{subequations}\label{eq:terms_wilson_limit}
\begin{align}
 \lim_{q^{a}\to\infty}w_x(a,b,c,d,e,f)&=q^{(1+\beta)x}
  \frac{(q^{-N};q)_x(q^{1-\alpha};q)_x}{(q;q)_x(q^{\beta-\alpha-N+2};q)_x}
  =\bar{w}_x(\alpha,\beta,N;q),\\
 \lim_{q^{a}\to\infty}u_n(x+a;a,b,c,d,e,f)
  &=\pPq{3}{2}{q^{-n},q^{n+\beta-N},q^{-x}}{q^{-N},q^{\alpha-x}}{q}{q^{\alpha-\beta}}
  =\bar{u}_n(x;\alpha,\beta,N;q),\\
 \lim_{q^{a}\to\infty}v_n(x+a;a,b,c,d,e,f)
  &=\pPq{3}{2}{q^{-n},q^{n+\beta-N},q^{x-N}}{q^{-N},q^{x-N+\beta-\alpha+2}}{q}{q}
  =\bar{v}_n(x;\alpha,\beta,N;q),\\
 \lim_{q^{a}\to\infty}h_n(a,b,c,d,e,f)
  &=q^{N(\alpha-1-n)}(q;q)_n\frac{(q^{-\beta};q)_N}{(q^{\alpha-\beta-1};q)_N}
  \frac{(q^{\beta+1};q)_n}{(q^{-N};q)_n}\frac{(q^{n+\beta-N};q)_n}{(q^{1+\beta-N};q)_{2n}}
  \no\\
  &=\bar{h}_n(\alpha,\beta,N;q)
\end{align}
\end{subequations}
and the biorthogonality relation is written
\begin{align}\label{eq:biortho_wilson_limit}
 \sum_{x=0}^{N}\bar{w}_x(\alpha,\beta,N;q)\, \bar{u}_n(x;\alpha,\beta,N;q)\,
 \bar{v}_n(x;\alpha,\beta,N;q) = \delta_{nm}\,\bar{h}_n(\alpha,\beta,N;q).
\end{align}
Note that
\begin{align}\label{}
 \bar{v}_n(x;\alpha,\beta,N;q)=\bar{u}_n(N-x;\beta-\alpha+2,\beta,N;q^{-1}).
\end{align}
Using \eqref{eq:terms_wilson_limit}-\eqref{eq:biortho_wilson_limit} and comparing with
\eqref{eq:weight}, \eqref{eq:Unx}, \eqref{eq:Vnx}, one can now extract $H_n$ in the
biorthogonality relation \eqref{eq:bio}.

\subsection{The $q\to1$ limit}\label{sec:qto1}
The $q=1$ limiting case of the present paper was previously investigated in
\cite{TsujimotoVinetetal2020}.  There, the authors presented the biorthogonality relation
of the BRF of Hahn type\footnote{Note that the proof of the biorthogonality given in
\cite{TsujimotoVinetetal2020} is not correct although the end result is right. The proof
in Section \ref{sec:ortho} of the present paper can be adapted to the $q\to1$ case.} but
they did not include the normalization.
We present it here.
It can be obtained straightforwardly as a $q\to1$ limit of
\eqref{eq:terms_wilson_limit}-\eqref{eq:biortho_wilson_limit}.
To that end, write $q=e^{h}$ and take the limit $h\to0$.
This leads to the biorthogonality relation
\begin{align}\label{eq:biortho_qto1_limit}
 \sum_{x=0}^{N}\underline{w}_x(\alpha,\beta,N)\, \underline{u}_n(x;\alpha,\beta,N)\,
 \underline{v}_n(x;\alpha,\beta,N) = \delta_{nm}\,\underline{h}_n(\alpha,\beta,N)
\end{align}
with
\begin{subequations}\label{eq:terms_qto1_limit}
\begin{align}
 \lim_{h\to0}\bar{w}_x(\alpha,\beta,N;e^h)&=
  \frac{(-N)_x(1-\alpha)_x}{x!(\beta-\alpha-N+2)_x}
  =\underline{w}_x(\alpha,\beta,N),\\
 \lim_{h\to0}\bar{u}_n(x;\alpha,\beta,N;e^h)
  &=\pFq{3}{2}{-n,n+\beta-N,-x}{-N,\alpha-x}{1}
  =\underline{u}_n(x;\alpha,\beta,N),\\
 \lim_{h\to0}\bar{v}_n(x;\alpha,\beta,N;e^h)
  &=\pFq{3}{2}{-n,n+\beta-N,x-N}{-N,x-N+\beta-\alpha+2}{1}
  =\underline{v}_n(x;\alpha,\beta,N),\\
 \lim_{h\to0}\bar{h}_n(\alpha,\beta,N;e^h)
  &=n!\frac{(-\beta)_N}{(\alpha-\beta-1)_N}
  \frac{(\beta+1)_n}{(-N)_n}\frac{(n+\beta-N)_n}{(1+\beta-N)_{2n}}
  =\underline{h}_n(\alpha,\beta,N).
\end{align}
\end{subequations}
\begin{rmk}
We note that Wilson also provided $q=1$ analogues of his ${}_{10}\phi_9$ BRF in
\cite{Wilson1991}, expressed as truncating hypergeometric series of ${}_9F_8$ type. It is
also possible to recover \eqref{eq:biortho_qto1_limit}-\eqref{eq:terms_qto1_limit}
starting instead from these functions, taking the parametrization in
\eqref{eq:parametrization} and then taking the limit $a\to\infty$.
\end{rmk}

\section{Tridiagonal actions of the triplet $X$, $Y$, $Z$}\label{sec:tridiag}
We have already shown how the operators $X$, $Y$, $Z$ act as matrices on the
basis $e_n(x)$ defined on the finite uniform grid $x=0,1,\dots,N$. The operator $Y$ is
tridiagonal, while $X$ and $Z$ are lower bi-diagonal (the upper off-diagonal is zero).

Alternatively, in the basis $\qphi{n}{x}{\alpha}$ the operator $Y$ is again a tridiagonal
matrix while $X$ and $Z$ are upper bi-diagonal.

Now consider the action of these operators on the rational functions
$U_n(x;\alpha,\beta,N;q)$. The functions $U_n$ form another basis on the
finite-dimensional space $M_N$. The actions of $X$, $Y$, $Z$ on this basis are given as
follows:
\begin{subequations}\label{eq:triplettridiag}
\begin{align}\label{}
 X\,U_n(x;\alpha,\beta,N;q)&=\mu_n^{(1)}U_{n+1}(x;\alpha,\beta,N;q)
  +\mu_n^{(2)}U_{n}(x;\alpha,\beta,N;q)+\mu_n^{(3)}U_{n-1}(x;\alpha,\beta,N;q),\\
 Y\,U_n(x;\alpha,\beta,N;q)&=\mu_n^{(4)}U_{n+1}(x;\alpha,\beta,N;q)
 +\mu_n^{(5)}U_{n}(x;\alpha,\beta,N;q)+\mu_n^{(6)}U_{n-1}(x;\alpha,\beta,N;q),\\
 Z\,U_n(x;\alpha,\beta,N;q)&=\mu_n^{(7)}U_{n+1}(x;\alpha,\beta,N;q)
  +\mu_n^{(8)}U_{n}(x;\alpha,\beta,N;q)+\mu_n^{(9)}U_{n-1}(x;\alpha,\beta,N;q),
\end{align}
where
\begin{align}
\begin{aligned}
 \mu_n^{(1)}&=-q^{-\alpha-n}
  \frac{[n]_q[\beta+n+1]_q[N-\beta-n]_q}{[N-\beta-2n]_q[2n+1+\beta-N]_q},\\
 \mu_n^{(2)}&=q^{-\alpha}\left(-[\alpha]_q+[-n]_q
  +\frac{[n]_q[1-n]_q[\beta+n]_q}{[2n-1+\beta-N]_q}
  -\frac{[-n]_q[n+1]_q[\beta+n+1]_q}{[2n+1+\beta-N]_q}\right),\\
 \mu_n^{(3)}&=-q^{-\alpha}
  \frac{[-n]_q[N-\beta-n]_q[N-n+1]_q}{[N-\beta-2n]_q[N-\beta-2n+1]_q},\\
 \mu_n^{(7)}&=q^{-\alpha-n}
  \frac{[\beta+n+1]_q[N-\beta-n]_q}{[N-\beta-2n]_q[2n+1+\beta-N]_q},\\
 \mu_n^{(8)}&
  =q^{-\alpha+1}\left(\frac{[-n-1]_q[N-n]_q}{[2n+1+\beta-N]_q}
   -\frac{[-n]_q[N-n+1]_q}{[2n-1+\beta-N]_q}\right)\\
  &\quad+q^{-\alpha-n}\left((1-q^{N+1})-q^{N-n}(1-q)\right)-1,\\
 \mu_n^{(9)}&=q^{-\alpha}\frac{[-n]_q[N-n+1]_q}{[N-\beta-2n]_q[N-\beta-2n+1]_q}.
\end{aligned}
\end{align}
As for the coefficients $\mu_n^{(4)}$, $\mu_n^{(5)}$, $\mu_n^{(6)}$, they are proportional
to the coefficients $\mu_n^{(1)}$, $\mu_n^{(2)}$, $\mu_n^{(3)}$
\begin{align}\label{}
 \mu_n^{(3+i)}=[-n]_q[n+\beta-N]_q\,\mu_n^{(i)},\qquad i\in \{1,2,3\}.
\end{align}
The formulas above are valid for $n=0,\dots,N$ except for $\ell\in\{1,4,7\}$, for which
the expressions for $n=N$ are given by
\begin{align}\label{}
 \mu_N^{(1)}=\mu_N^{(4)}=\mu_N^{(7)}=0.
\end{align}
\end{subequations}

\section{Bispectrality}\label{sec:bispectrality}

Collecting our observations, let us now record how the bispectrality of the BRF of
$q$-Hahn type is expressed through two generalized eigenvalue equations: a difference
equation in the variable $x$ with eigenvalue depending on the degree $n$, and a recurrence
relation in the variable $n$ with eigenvalue depending on $x$.

\subsection{The difference equation}\label{sec:diffeqn}
Return to the GEVP \eqref{eq:gevp} that the operators $X$ and $Y$ define.
We have established in Proposition \ref{thm:gevp} that the $U_{n}(x;\alpha,\beta,N;q)$
are its solutions, with eigenvalues ${\lambda_n=[-n]_q[n+\beta-N]_q}$.

From the expressions of $X$ and $Y$ given in \eqref{eq:Xdef} and
\eqref{eq:Ydef}-\eqref{eq:coeffYshift} we see that this amounts to the equation
\begin{align}\label{}
\begin{aligned}
 &A_1(x)U_{n}(x+1;\alpha,\beta,N;q)+A_0(x)U_{n}(x;\alpha,\beta,N;q)+
 A_2(x)U_{n}(x-1;\alpha,\beta,N;q)\\
 &\quad=\lambda_n\left(
 [x-\alpha]_qU_{n}(x;\alpha,\beta,N;q)
 -q^{-\alpha}[x]_qU_{n}(x-1;\alpha,\beta,N;q)\right),
\end{aligned}
\end{align}
with $A_0(x)$, $A_1(x)$, $A_2(x)$ given in \eqref{eq:coeffYshift}.
This is the difference equation for the $U_{n}(x;\alpha,\beta,N;q)$ BRF.

\subsection{The recurrence relation}\label{sec:recrel}
Next, we look at the pair of operators $X$ and $Z$.
As per \eqref{eq:zfx}, these define the GEVP
\begin{align}\label{eq:gevpxz}
 X\,U_n(x)=-[x-\alpha]_qZ\,U_n(x).
\end{align}
From the
actions on the $U_n$ basis given in the preceding section we see that \eqref{eq:gevpxz}
amounts to the recurrence relation
\begin{align}\label{}
 &\mu_n^{(1)}U_{n+1}(x;\alpha,\beta,N;q)+\mu_n^{(2)}U_{n}(x;\alpha,\beta,N;q)+
 \mu_n^{(3)}U_{n-1}(x;\alpha,\beta,N;q)\\
 &\quad=-[x-\alpha]_q\left(
 \mu_n^{(7)}U_{n+1}(x;\alpha,\beta,N;q) +\mu_n^{(8)}U_{n}(x;\alpha,\beta,N;q)
 +\mu_n^{(9)}U_{n-1}(x;\alpha,\beta,N;q)\right)\no
\end{align}
for the BRF $U_{n}(x;\alpha,\beta,N;q)$.

\section{An algebra for the bispectrality of the $q$-Hahn rational functions}
\label{sec:algebra}
In the previous section, we have used the three operators $X$, $Y$, $Z$ to obtain the
bispectral equations of the rational functions.  For orthogonal polynomials of the
Askey scheme, one would use two operators to encode the bispectrality through
eigenvalue equations. The algebraic relations obeyed by these two operators lead to rich
algebraic structures. This is precisely how one obtains the Askey--Wilson algebra
relations, starting from the operators defining the difference and recurrence relations of
the Askey--Wilson polynomials \cite{Zhedanov1991}.  To extend this
construction, we now look at the algebraic relations obeyed by the triplet of operators
that express the bispectrality properties of the biorthogonal rational
functions $U_n$, the operators $X$, $Y$, $Z$.

The relations realized by $X$, $Y$, $Z$ will be taken to represent the following
abstract algebra.

\begin{defi}
The rational $q$-Hahn algebra $r\mathfrak{h}_q$ is the unital associative algebra
generated by $X$, $Y$, $Z$ obeying the relations
\begin{align}
 XZ-qZX&=-Z^2-Z+\xi_6X,\label{eq:comXZ}\\
 ZY-qYZ&=\xi_1X^2+\xi_3\{X,Z\}+\xi_4Z^2+\xi_5X+\xi_6Y+\xi_7Z+\xi_0,\label{eq:comZY}\\
 YX-qXY&=\xi_3X^2+\xi_4\{X,Z\}-\{Y,Z\}+\xi_2Z^2+\xi_7X-Y+\xi_8Z+\xi_0,\label{eq:comYX}
\end{align}
with
\begin{gather}\label{}
\begin{gathered}
 \xi_1=-q^{\beta-N-1}[3]_q,\qquad\xi_2=q^{-1}[\beta-N]_q,\\
 \xi_3=[\beta-N-1]_q+q^{\beta-N},\qquad
 \xi_4=-q^{\beta-N-1},\\
 \xi_5=[\beta-\alpha]_q+[\beta-\alpha+1]_q
  +q^{\beta-N}\left(1+[-\alpha]_q+[-\alpha-1]_q\right),\\
 \xi_6=-(1-q),\qquad\xi_7=q^{\beta-N}[-\alpha-1]_q
  -[-\alpha]_q[\beta-N]_q,\\
 \xi_8=q^{-1}[\beta-\alpha-N]_q
  -[-\alpha]_q[\beta-N]_q,\qquad
 \xi_0=-[-\alpha]_q[\beta-\alpha]_q.
\end{gathered}
\end{gather}
\end{defi}
\noindent
It can be observed that this algebra derives from a potential $\Phi_{r\mathfrak{h}_q}$
\cite{IyuduShkarin2018}.
Let $F=\mathbb{C}[x_1,\dots,x_n]$ be a free associative algebra with $n$
generators and $F_{cyc}=F/[F,F]$. The basis of $F_{cyc}$ is spanned by the cyclic words
$[x_{i_1}\dots x_{i_k}]$. One can introduce the map $\tfrac{\partial}{\partial{x_j}}:
F_{cyc}\to F$ that acts as follows on basis elements:
\begin{align}\label{}
 \frac{\partial[x_{i_1}\dots x_{i_r}]}{\partial{x_j}}
 =\sum_{s|i_s=j}x_{i_{s+1}}x_{i_{s+2}}\dots x_{i_r}x_{i_1}x_{i_2}\dots x_{i_{s-1}}.
\end{align}
An algebra $F$ whose defining relations are given by
\begin{align}\label{}
 \frac{\partial\Phi}{\partial{x_j}}=0,\qquad j=1,\dots,n
\end{align}
for a $\Phi\in F_{cyc}$ is said to derive from a potential $\Phi$.

In the case of the rational $q$-Hahn algebra $r\mathfrak{h}_q$, it is seen that the
defining relations are obtained from the following potential:
\begin{align}
\begin{aligned}
 \Phi_{r\mathfrak{h}_q}&=q[XYZ]-[YXZ]+\frac{\xi_1}{3}[X^3]+\frac{\xi_2}{3}[Z^3]
  +\xi_3[X^2Z]+\xi_4[XZ^2]-[YZ^2]\\
  &\quad +\frac{\xi_5}{2}[X^2]+\xi_6[XY]+\xi_7[XZ]-[YZ]+\frac{\xi_8}{2}[Z^2]
  +\xi_0\left([X]+[Z]\right).
\end{aligned}
\end{align}
More precisely, \eqref{eq:comXZ}-\eqref{eq:comYX} are respectively equivalent to
\begin{align}\label{}
 \frac{\partial\Phi_{r\mathfrak{h}_q}}{\partial{Y}}=0,\qquad
 \frac{\partial\Phi_{r\mathfrak{h}_q}}{\partial{X}}=0,\qquad
 \frac{\partial\Phi_{r\mathfrak{h}_q}}{\partial{Z}}=0.
\end{align}
A Casimir element for this algebra is given by
\begin{align}
\begin{aligned}
 Q_{r\mathfrak{h}_q}&=(1-q)XYZ+q^{\beta-N+1}X^3+\gamma_{1}X^2Z+\gamma_{2}XZ^2+qYZ^2\\
  &\quad+\gamma_{3}X^2+(1-q)XY+\gamma_{4}XZ-(1-2q)YZ+\gamma_{5}X-(1-q)Y
\end{aligned}
\end{align}
with
\begin{gather}\label{}
\begin{gathered}
 \gamma_{1}=-q^{\beta-N}-[\beta-N+2]_q,\qquad
 \gamma_{2}=[\beta-N+1]_q,\\
 \gamma_{3}=-q^{\beta-N}[-\alpha]_q-[\beta-\alpha+1]_q
  -q^{\beta-N+1},\\
 \gamma_{4}=q^{\beta-N}+[\beta-\alpha]_q+[\beta-\alpha+1]_q
  +[-\alpha]_q[\beta-N+1]_q,\\
 \gamma_{5}=[\beta-\alpha+1]_q+q^{\beta-N}[-\alpha]_q
  +[-\alpha]_q[\beta-\alpha]_q.
\end{gathered}
\end{gather}
Next, let us bring to the fore the factorization of $Y$ as $XV$.
In the $q$-shift representation, the operator $V$ takes the form
\begin{align}
\begin{aligned}
 V&=q^{-x+\beta}[x-N]_q[x-\alpha+1]_qT^{+}\\
  &\quad +q^{-\alpha+\beta+1}[x]_q[x-N-1]-q^{-x+\beta}[x-N]_q[x-\alpha+1]_q
  -q^{-x}[x]_q[x-N-\alpha+\beta]_q \\
  &\quad +q^{\beta+1-\alpha-x}[\alpha-\beta-1]_q[\alpha-1]_q
  \sum_{k=1}^{N}q^{k}\frac{(q^{-x};q)_k}{(q^{\alpha-x};q)_k}\left(T^{-}\right)^k,
\end{aligned}
\end{align}
with $(T^{-})^{N+1}$ acting as zero in our truncated grid.
Based on the nomenclature introduced in the $q=1$ case \cite{VinetZhedanov2020a}, we
introduce what we will call the meta-$q$-Hahn algebra $m\mathfrak{h}_q$.
\begin{defi}
The meta-$q$-Hahn algebra $m\mathfrak{h}_q$ is the unital associative algebra
generated by $X$, $V$, $Z$ obeying the relations
\begin{align}
 XZ-qZX&=-Z^2-Z-(1-q)X,\\
 VX-qXV&=-\{V,Z\}+\eta_1X-V-q^{-1}\eta_1Z-\eta_0,\\
 ZV-qVZ&=-\eta_2X-(1-q)V+\eta_1Z+\eta_3,
\end{align}
with
\begin{align}
\begin{aligned}
 \eta_0&=[-\alpha]_q[\beta-N]_q+q^{-1}[\beta-\alpha+1]_q,\\
 \eta_1&=[\beta-N+1]_q,\\
 \eta_2&=-q^{\beta-N}[2]_q,\\
 \eta_3&=[\beta-\alpha+1]_q+q^{\beta-N}[-\alpha]_q.
\end{aligned}
\end{align}
\end{defi}
\noindent
We can compute a Casimir element for this algebra which is cubic:
\begin{align}\label{}
\begin{aligned}
 Q_{m\mathfrak{h}_q}&=(1-q)XVZ+qVZ^{2}+ q^{\beta-N+1}X^{2}+(1-q)XV+\gamma_1XZ
 -(1-2q)VZ\\
 &\quad+\eta_1Z^{2}+\gamma_3X-(1-q)V+\gamma_4Z.
\end{aligned}
\end{align}
One can see that the $q$-commutators between $V$ and other generators are linear on the
rhs up to a single quadratic term; they thus depart minimally from pure linear
relations.
Such a characterization, and the fact that the algebra is connected to rational functions
of $q$-Hahn type, hints that the meta-$q$-Hahn algebra is bound to have interesting
properties.  As an example, let us mention that it derives from the potential
$\Phi_{m\mathfrak{h}_q}$:
\begin{align}
\begin{aligned}
 \Phi_{m\mathfrak{h}_q}&=q[XVZ]-[VXZ]-[VZ^2]-(1-q)[XV]-[VZ]+\eta_1[XZ]
 -\frac{q^{-1}\eta_1}{2}[Z^2]-\frac{\eta_2}{2}[X^2]\\
 &\quad+\eta_0[Z]+\eta_3[X].
\end{aligned}
\end{align}
Its representation theory will then be intimately tied to the rational functions of
$q$-Hahn type.

\section{Conclusion}\label{sec:conclusion}
Let us repeat the main points of this paper.  We have introduced a triplet of
$q$-difference operators $X$, $Y$, $Z$ that play a role analoguous to the bispectral
operators of orthogonal polynomials in the case of biorthogonal rational functions. We
have shown how these three operators provide a characterization of the bispectrality
properties of the BRF of $q$-Hahn type through GEVP, and we have also proven the
biorthogonality of these functions.  The underlying algebra
was obtained and the connection with Wilson's ${}_{10}\phi_9$ BRF was made explicit. A
remarkable factorization of the operator $Y$ as $XV$ was exhibited.

In future work we plan to explore the representations of the meta-$q$-Hahn algebra
introduced in the last section.  It is expected that this algebra will provide a unified
algebraic interpretation of the bispectrality properties of both the $q$-Hahn orthogonal
polynomials and the biorthogonal rational functions of $q$-Hahn type.

This work adds to the systematic exploration and characterization of the bispectrality
of BRF. The ${}_{3}\phi_2$ and ${}_{3}F_2$ Hahn type BRF have now been well studied
\cite{VinetZhedanov2019, TsujimotoVinetetal2020, VinetZhedanov2020a} and
the next step would be to explore BRF of $q$-Racah type. These should be expressable as
${}_{4}\phi_3$ and their $q\to1$ limits should provide BRF of Racah type. Note that the
Askey biorthogonal polynomials on the unit circle have also been shown in
\cite{VinetZhedanov2021} to be amenable to a similar treatment.

Lastly, let us mention that Wilson's BRF have appeared in the study of the so-called
rational Heun operators \cite{TsujimotoVinetetal2019}.  These rational Heun operators were
defined by extending the construction of Heun operators to rational functions.
Sklyanin--Heun operators have been recently introduced as building blocks of the Heun
operators and have led to an algebraic characterization of various families of
para-polynomials \cite{GaboriaudTsujimotoetal2020, BergeronGaboriaudetal2021,
BergeronGaboriaudetal2021a}. One can expect a similar story to unfold for the rational
functions.  Rational analogues of Sklyanin--Heun operators should prove interesting
objects of study and should lead to the definition of para families of rational functions.

\subsection*{Acknowledgments}\label{sec:ack}
I. Bussière benefitted from a CRM-ISM intern scholarship.
J. Gaboriaud received financial support from the ESP of the Université de Montréal.
The research of L. Vinet is supported in part by a Discovery Grant from NSERC.
A. Zhedanov who is funded by the National Foundation of China (Grant No.11771015)
gratefully acknowledges the hospitality of the CRM over an extended period.

\printbibliography

\end{document}